\documentclass[12pt]{article}
\usepackage{dotf,float,amssymb,geompsfi}

\begin{document}
\def\citeauthoryear#1#2#3{#1, #3}

\date{}
\title{Akbulut's corks and h-cobordisms of smooth, simply connected
4-manifolds}
\author{Rob Kirby}
\maketitle


Theorem:  Let $M^5$ be a smooth 5-dimensional h-cobordism between two
simply connected, closed 4-manifolds, $M_0$ and $M_1$.  Then there
exists a sub-h-cobordism $A^5 \subset M^5$ between $A_0 \subset M_0$
and $A_1 \subset M_1$ with the properties:

\begin{description}

\item[(1)] $A_0$ and hence $A$ and $A_1$ are compact contractible
manifolds, and

\item[(2)] $M-intA$ is a product h-cobordism, i.e. it is
diffeomorphic to $(M_0-intA_0) \times [0,1]$.

\end{description}

This theorem first appeared in a preprint of Curtis \& Hsiang in fall
1994.  Soon after, much shorter proofs were found by Freedman \& Stong
\cite{cur95}, Matveyev \cite{mat95}, and Z. Bi\v{z}aca.  The following
improvements were also shown:

Addenda: The h-cobordism $A$ can be chosen so that,

\begin{description}

\item[(A)] $M-A$ (and hence each $M_i - A_i$) is simply connected
(Freedman \& Stong) \cite{cur95},

\item[(B)] $A$ is diffeomorphic to $B^5$ (Bi\v{z}aca, Kirby) (but not, of
course, preserving the structure of the h-cobordism),

\item[(C)] $A_0 \times I$ and $A_1 \times I$ are diffeomorphic to $B^5$
\cite{mat95},

\item[(D)] $A_0$ is diffeomorphic to $A_1$ by a diffeomorphism which,
restricted to $\partial A_0 = \partial A_1$, is an involution
\cite{mat95}.

\end{description}

Corollary:  Any homotopy 4-sphere, $\Sigma^4$, can be constructed by
cutting out a contractible 4-manifold, $A_0$ from $S^4$ and gluing it
back in by an involution of $\partial A_0$.

Remark: Since there are many examples of non-trivial h-cobordisms (the
first ones were discovered by Donaldson \cite{don87}), there are as
many examples of non-trivial, rel boundary, h-cobordisms $A$.  However
these  $A$ are delicate objects;  their non-triviality vanishes when a
trivial h-cobordism is added.  That is, if we add $A_0 \times I$ to
$A$ along $\partial A_0 \times I$, then it follows from the Addenda
that we have an h-cobordism between $S^4$ on the bottom as well as
$S^4$ on the top;  thus the h-cobordism is the trivial  $S^4 \times
I$.

\begin{figure}[H]
\centerline{\psfig{figure=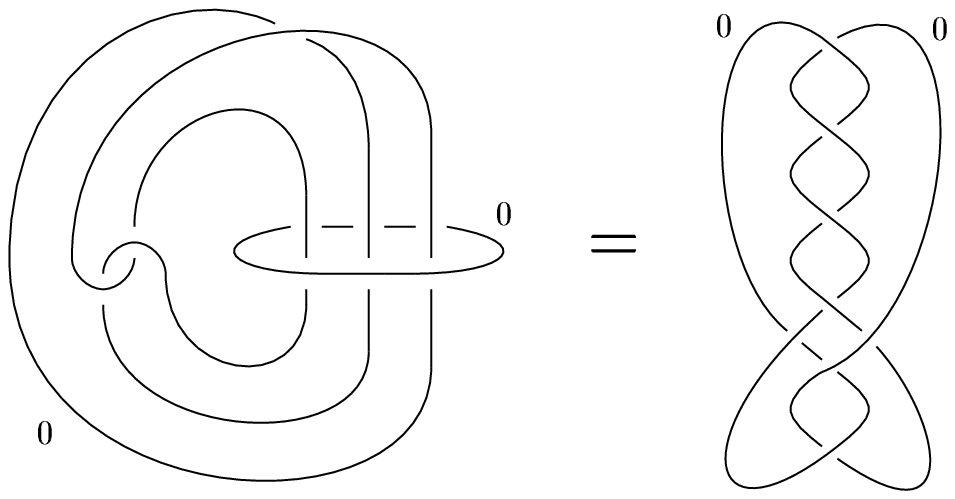,height=2.0in}}
\caption{}
\label{acork}
\end{figure}

{\bf 1.} The first example of a non-trivial h-cobordism, $A$, on $B^5$
was found by Akbulut \cite{akb91b}.  It is the prototype of the
h-cobordisms $A$ in the theorem, and it seems appropriate to call such
h-cobordisms {\it Akbulut's corks}, for any exotic h-cobordism can be
constructed from the product h-cobordism by pulling out a {\it cork}
and putting it back in with a {\it twist}, (which preserves $A_0$ but
not the structure of the h-cobordism).

Akbulut constructs a homology $S^2 \times S^2 - point$, called
$A_{1/2}$, by adding two 2-handles to the symmetric link $L$ of
unknots drawn in Figure~\ref{acork}.  Because this link $L$ is symmetric, there
is an involution $\tau : \partial A_{1/2} \rightarrow \partial
A_{1/2}$, which extends to a diffeomorphism $\overline{\tau} : A_{1/2}
\rightarrow A_{1/2}$ which switches the 2-handles (to extend $\tau$
over the 0-handle of $A_{1/2}$, one can either cone the involution on
$S^3$ obtaining an involution which is not smooth at the cone point,
or extend over the 0-handle using the fact that the involution on
$S^3$ is diffeotopic to the identity).

Since the components of $L$ are 0-framed unknots, one can trade either
of the 2-handles (but not both for $L$ is not the unlink) for a
1-handle, obtaining $A_0$ or $A_1$ respectively. (Recall that adding a
2-handle to a 0-framed unknot gives $S^2 \times B^2$, adding an
orientable 1-handle to $B^4$ gives $B^3 \times S^1$, and both have
boundary $S^2 \times S^1$; thus one may change the 4-manifold by trading
2-handles for 1-handles, or vice versa, and then adding all other
handles to the $S^2 \times S^1$ boundary as before the trade.)  The
new 1-handles are denoted by the same unknot but with a dot on it,
which means that any arcs going through the dotted circle actually go
over a 1-handle.  This can be seen by observing that a dotted circle
means: remove the obvious, properly imbedded, 2-ball in the 0-handle
leaving $S^1 \times B^3$, which is a 0-handle union a 1-handle.  This
operation does not change boundaries, so $\partial A_{1/2} = \partial
A_0 = \partial A_1$.

Akbulut proves (using a long series of handle moves
culminating in an application of Donaldson's invariants) that the
identity map, $id:\partial A_0 \rightarrow \partial A_1$, does not
extend to a diffeomorphism of $A_0$ to $A_1$.

The notation $A_0$, $A_{1/2}$, $A_1$, suggests there is an h-cobordism
lurking about, and this is correct.  The operation of trading
2-handles for 1-handles can also be done by adding 3-handles in the
right way.  Each component of $L$, $K_0$ or $K_1$, determines a
2-sphere, $S_0$ or $S_1$, composed of the core of the 2-handle and the
obvious slice disk that $K_i$ bounds in the 0-handle.  Each 2-sphere
has a trivial normal bundle (because of the 0-framing), and $S_0 \cap
S_1$ is three points (algebraically one).  To construct an h-cobordism
$A$, start with $A_{1/2} \times [1/2 -\epsilon, 1/2 + \epsilon]$ and
add a 3-handle to $S_0 \times (1/2 - \epsilon)$ and a 3-handle to
$S_1 \times (1/2 + \epsilon)$.  The new boundary on the bottom will be
$A_0$ because $S_0 \times B^2$ has been removed from $A_{1/2} \times
(1/2 - \epsilon)$, thus removing the 2-handle and the slice disk
which has the effect of switching the 2-handle to a 1-handle.
Similarly for $A_1$.  The structure of the h-cobordism $A$ is to add a
2-handle to $A_0$ (the 3-handle turned upside down) and the 3-handle
to $S_1$.  $A$ must be non-trivial because it is  a product over
$\partial A_0$, describing $id: \partial A_0 \rightarrow \partial
A_1$, and Akbulut showed this cannot extend to a diffeomorphism from
$A_0$ to $A_1$.

There is a natural generalization of Akbulut's cork.  Let the 0-handle
in $A_{1/2}$ be replaced by a contractible 4-manifold $B_{1/2}$.
Suppose $D$ is a collection of $2n$ properly imbedded 2-balls,
$D_{0,i} \cup D_{1,i}, i = 1\ldots n$, in $B_{1/2}$, and suppose that
$D_{0,i} \cap D_{0,j} = \emptyset = D_{1,i} \cup D_{1,j}$ for all $i,j
\in {1, \ldots ,n}$, and that algebraically $D_{0,i} \cap D_{1,j} =
\delta_{ij}$.  Then we can form $A_{1/2}$ by adding 2-handles with
0-framings to each $\partial D_{0,i}$ and $\partial D_{1,i}$, $i = 1,
\ldots , n$.  This produces obvious 2-spheres with trivial normal
bundles, $S_{0,i}$ and $S_{1,i}$, $i = 1, \ldots, n$.  Then we form
our h-cobordism $A$ by adding $n$ 3-handles below and $n$ 3-handles
above to the $S_{0,i}$'s and the $S_{1,i}$'s respectively.

We can conjecture that the product structure on the sides of $A$,
namely $\partial A_0 \times I$, extends over $A$ iff $D$ is
{\it concordant} in $B_{1/2} \times I$ to the $\partial$-connected sum
of $n$ copies of $B^2 \times 0 \cup 0 \times B^2$ in $B^4 \subset B_{1/2}$.

{\bf 2.} Here is a proof of the Theorem.  The exposition is not
particularly original, but gains by organizing all the steps into a
whole rather than having them split into the two papers \cite{cur95}
and \cite{mat95}.

We begin with a Morse function $f:(M, M_0, M_1) \rightarrow (I,0,1)$
and its associated handlebody structure which adds k-handles, $0
\leq k \leq 5$, to $M_0$.  We can cancel all 0-handles and 5-handles
since $M$ is connected.  We can cancel all 1-handles and 4-handles (at
the cost of new 2 and 3-handles) just as Smale did in the original
proof of the higher dimensional h-cobordism theorem (\cite{rou72} Lemma 6.15).
Alternatively, we may have the h-cobordism provided by Wall (see
\cite{wal64b}, \cite{kir89} Chapter 9) which begins with homotopy
equivalent, simply connected closed, smooth 4-manifolds $M_0$ and
$M_1$ and constructs $M$ using only 2 and 3-handles.  (This involves
no loss of generality because any two h-cobordisms between $M_0$ and
$M_1$ are diffeomorphic \cite{kre94,kre95}.)

We can assume that $M_{1/2} = f^{-1}(1/2)$ has all the 2-handles below
and all the 3-handles above.  Note that $M_{1/2}$ is 1-connected which
is always true of the upper boundary when 2-handles are attached to a
\{simply connected 4-manifold\} $\times I$.  Each 2-handle has an
ascending 3-ball which meets $M_{1/2}$ in a smoothly imbedded 2-sphere;
call these $S_{0,i}, i = 1, \ldots , n$.  Similarly each 3-handle
descends to meet $M_{1/2}$ in $S_{1,i}, i = 1, \ldots , n$.  We can
assume, perhaps after some handle slides, that the boundary map from
3-chains (generated by the 3-handles) to 2-chains (generated by
2-handles) is given by the identity matrix, or, equivalently, that
algebraically $S_{0,i} \cap S_{1,j} = \delta_{ij}$.

\begin{figure}
\centerline{\psfig{figure=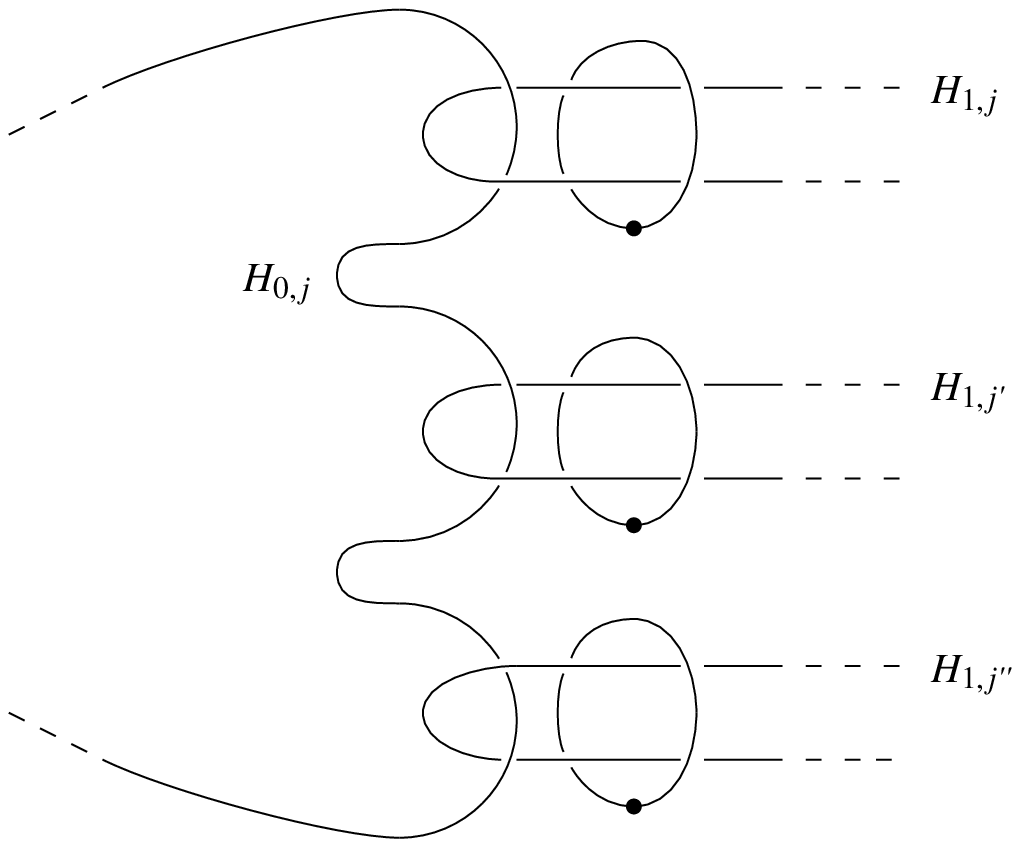}}
\caption{}
\label{1-handles}
\end{figure}

{\bf 3.} Choose a base point $\ast$ in $M_{1/2}$ minus all the spheres
$S_{k,i}$.  Choose $2n$ arcs in general position which connect $\ast$
to basepoints $\ast_{k,i}$ in the $S_{k,i}$.  A regular neighborhood
of these arcs will be our 0-handle in a forthcoming handlebody
structure on $M_{1/2}$.  In each $S_{k,i}$, run disjoint arcs from
$\ast_{k,i}$ to each point of intersection of $S_{k,i}$ with some
$S_{k^{\prime},j}$ (with $k \neq k^{\prime}$).  These arcs come in
pairs, from $\ast_{k,i}$ to a point in $S_{k,i} \cap S_{k^{\prime},j}$
to $\ast_{k^{\prime},j}$, and a regular neighborhood forms a 1-handle
attached to the 0-handle.  Each $S_{k,i}$, minus the regular
neighborhood of the tree of arcs in it, gives a 2-handle, $H_{k,i}$,
which is added to the 0-handle and 1-handles.  Figure~\ref{1-handles} shows 
how the $H_{k,i}$'s behave with respect to the 1-handles; note that each
$H_{k,i}$ is attached to an unknot and the $H_{0,i}$'s and the
$H_{1,i}$'s are each attached to unlinks of $n$ components.  Note that
we can assume that the $H_{0,i}$ do not go over any 1-handles, whereas
the $H_{1,i}$ go over and back so that, if the one handles correspond
to generators $x_1, x_2, \ldots x_r$ of the fundamental group $\pi_1
(M_{1/2})$, then the $H_{1,i}$ give relators equal to a product of
$x_i\overline{x}_i$'s and $\overline{x}_i x_i$'s, $i \in {1, \ldots ,
r}$.

So far we have chosen a 0-handle and some 1 and 2-handles in
$M_{1/2}$.  Extend this handlebody to a handlebody structure on all of
$M_{1/2}$ (it may have more 1-handles (still indexed by $1\ldots r$)
as well as 2 and 3-handles, but extra 0- and 4-handles may be avoided).
Since $M_{1/2}$ is 1-connected and the $H_{k,i}$ give trivial
relators, it follows that the other 2-handles $H_l, l = 1, \ldots s$
must homotopically kill the 1-handles, where we assume that the
attaching circle of each $H_l$ has a base point $\ast_l$ which has
been connected by an arc to $\ast$.  Note that when we slide $H_l$
over $H_m$, along an arc $\lambda$ joining $\ast_l$ to $\ast_m$ then
we replace the relator $r_l$ by the relator $r_l \lambda r_m^{\pm 1}
\overline{\lambda}$; $\lambda$ can be chosen to be trivial if
necessary.

It follows from elementary combinatorial group theory that we can
slide 2-handles over 2-handles so as to end up with the $H_l, l =
1, \ldots r$, exactly killing the generators $x_1, \dots x_r$;
that is, $H_i = x_i w_i$ where $w_i$ cancels away to $1$ using
only the relations $x_j \overline{x}_j = 1 = \overline{x}_j x_j$.
During this process, it may have been necessary to add cancelling
pairs of 2- and 3-handles, so as to slide a new 2-handle over some
$H_l$ which is about to be altered by sliding over another handle; the
new 2-handle preserves the relator $r_l$ for later use in the sequence
of Tietze moves which reduces the original presentation to the {\it
trivial} one.  These new 2-3 pairs may be necessary to avoid the
difficulties inherent in the Andrews-Curtis Conjecture \cite{and65},
\cite{kir95} Problem 5.2.

Let $B_{1/2}$ be the contractible manifold formed by the 0-handle, all
the 1-handles, and the 2-handles $H_l, l= 1, \ldots , r$.  Let
$A_{1/2}$ be $B_{1/2}$ union the 2-handles $H_{k,i}, k \in {0,1}, i
\in {1,\ldots n}$.  Then $A$ will be $A_{1/2}$ (thickened by crossing
with $[1/2 - \epsilon, 1/2 + \epsilon]$) together with the 3-handles
added below to the $S_{0,i}$'s and above to the $S_{1,i}$'s.

Clearly $A$ is contractible (since $B_{1/2}$ is contractible and the
3-handles cancel the $H_{k,i}$).  Since $A$ contains the 2 and
3-handles of the h-cobordism $M$, it follows that $M - intA$ is a
product h-cobordism.  This finishes the proof of the Theorem.

{\bf 4.} Proof of the Addenda:

(B), (C) and (D) are easiest to prove so we start there. 

$A$ is diffeomorphic to $B_{1/2} \times I$ because the 3-handles added
to $A_{1/2}$ geometrically cancel the 2-handles $H_{k,i}$,
$k = 0,1$, $i = 1, \ldots n$.  Furthermore $B_{1/2} \times I$ is
diffeomorphic to $B^5$ because homotopic circles in a 4-manifold are
isotopic, so the attaching maps of the $H_l$'s can be isotoped to
geometrically cancel the 1-handles, leaving only the 0-handle of
$B_{1/2} \times I$. This proves (B).

$A_0$ is contractible because $A$ is, but we need to also know that
each 1-handle of $A_0$ is homotopically cancelled by a 2-handle. $A_0$
is $A_{1/2}$ but with a dot on each attaching circle of the
$H_{0,i}$'s.  These dotted circles give $n$ new generators, $y_i$, $i
= 1, \dots , n$, to the presentation for $\pi_1 (B_{1/2})$, and the
$H_{1,i}$, $i = 1, \ldots , n$, are $n$ new relators, $s_i$, $i = 1,
\ldots , n$.  At this point we need to go back and make a careful
choice of the arcs in each $S_{1,i}$ which join $\ast_{1,i}$ to the
points of intersection of $S_{1,i}$ with the spheres $S_{0,j}$, $j =
1, \ldots , n$.  We first run arcs from $\ast_{i,1}$ to all points of
intersection with $S_{0,1}$, then with $S_{0,2}$, then $S_{0,3}$, and
so on to $S_{0,n}$.  This is easy to do because trees do not separate
points in dimension 2.  With this choice of arcs, it follows that the
attaching circle of $H_{1,i}$ reads off the word $w_1 w_2 \ldots w_n$
where $w_j$ is a word in the $y_j$ and $\overline{y}_j$ with exponent
sum zero if $j \neq i$ and exponent sum one if $j = i$.

Thus the 2-handles $H_l$, $l = 1, \ldots , r$, kill $x_1, \ldots ,
x_r$ and then the 2-handles $H_{1,i}$, $i = 1, \ldots , n$, kill the
generators $y_1, \ldots , y_n$.  Therefore, $A_0 \times I$ is
diffeomorphic to $B^5$, because homotopy implies isotopy for
1-manifolds in 4-manifolds, so the 2-handles geometrically cancel the
1-handles since they do so homotopically.  Similarly $A_1 \times I $
is $B^5$.  This finishes the proof of Addenda (C).

{\bf 5.} To prove (D), we increase the size of the h-cobordism $A$.  Choose a
4-ball $B_0^4$ in $M_0$ such that $B^4_0 \cap A_0 = \partial B^4_0
\cap \partial A_0 = B^3$.  $M$ is a product, $B^4 \times I$, over
$B^4_0$. 

Since $A$ is $B^5$, it follows that $\partial A = A_0 \cup_{\partial}
A_1 = S^4$.  If we remove an open 4-ball, which intersects $\partial
A_i$ in a 3-ball, from $\partial A$, then the result, $(A_0
\cup_{\partial} A_1)_0$, can be identified with $B^4_0$.  Similarly,
using the fact that $A_0 \times I$ is $B^5$, we can identify $(A_0
\cup_{\partial} A_0)_1$ with $B^4_1$.  Then the product h-cobordism
$(B^4 \times I, B^4_0, B^4_1)$ can be identified with $(A_0 \times I
\cup_{\partial} A^{-1} , (A_0 \cup_{\partial} A_1)_0, (A_0
\cup_{\partial} A_0)_1 )$ where $A^{-1}$ is $A$ upside down and $A_0
\times I$ and $A^{-1}$ are joined along $(\partial A_0 -intB^3 )
\times I$.

\begin{figure}
\centerline{\psfig{figure=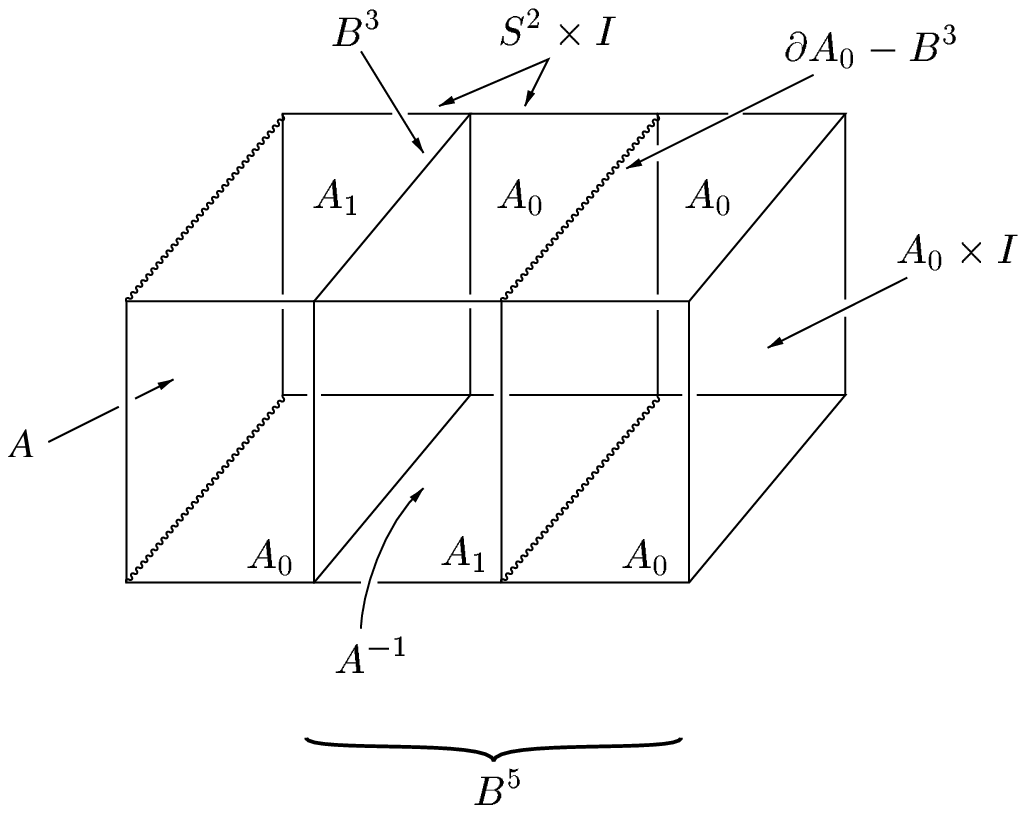}}
\caption{}
\label{symcork}
\end{figure}

Now we enlarge the h-cobordism $A$ by adding $A^{-1}$ to it (see
Figure~\ref{symcork}).  Clearly the complement is still a product, and clearly
the top and bottom of $A \cup A^{-1}$, namely $(A_1 \cup_{B^3} A_0)_1$
and $(A_0 \cup_{B^3} A_1)_0$, are diffeomorphic by the obvious
involution.  This proves (D).

{\bf 6.} To prove Addendum (A), that $A$ can be chosen so that the complement
$C = M - intA$ is simply connected, we must go back to the point in
the argument in which $A_{1/2}$ was constructed with {\it r}
1-handles, {\it r} 2-handles $H_l$, $l = 1, \ldots , r$, and the $2n$
2-handles $H_{k,i}, k \in {0,1}, i \in {1,\ldots n}$.  The complement
of $A_{1/2}$ has zero first homology, but it may not be simply
connected.

Let $L_1$ be a level set of $M_{1/2}$ after the 1-handles have been
attached to the 0-handle, ($L_1 = \# rS^1 \times S^2$), and let $L_3$
($= \# tS^1 \times S^2$) be a level set of $M_{1/2}$ just before the
3-handles are attached (equivalently, the boundary of the 4-handle
union the 3-handles).  Let $L_1 \cap L_3$ be denoted by $Q^3$; it can
be thought of as $L_1$ minus the attaching circles of all the
2-handles, or $L_3$ minus the co-circles of the 2-handles.  

Let $C_{1/2} = M_{1/2} - A_{1/2}$. Since $H_1 (C_{1/2}) = 0$, it
follows that $\pi_1 (C_{1/2})$ is generated by commutators, so we can
change it to zero if we have a method of sliding 2-handles that gives
us new 2-handles which kill commutators, but does not affect $\pi_1
(A_{1/2}) = 0$.  Here is such a method:

All slides of 2-handles over other 2-handles must take place along
arcs $\lambda$ lying in $Q$ (with endpoints at $\ast_{\alpha}$ and
$\ast_{\beta}$, which are connected to $\ast$ for fundamental group
computations).  If $H_{\alpha}$ and $H_{\beta}$ are 2-handles giving
relations $r_{\alpha}$ and $r_{\beta}$ in the generators $x_1 \ldots
x_r$ of $\pi_1 (L_1)$, and if we slide $H_{\alpha}$ over $H_{\beta}$
using the arc $\lambda$ and then slide $H_{\alpha}$ back over
$H_{\beta}$ using the arc $\mu$, then $r_{\alpha}$ is replaced by
$$r_{\alpha}\lambda r_{\beta} \lambda^{-1} \mu \overline{r}_{\beta}
\mu^{-1} = r_{\alpha} [\overline{\mu} \lambda, r_{\beta}]^{\mu}.$$
(Note that if $\lambda$ is homotopic to $ \mu$ in $\pi_1 (L_1)$, then
$r_{\alpha}$ is unchanged.)  The effect of these two slides on the
generators of $\pi_1 (L_3)$ is this: the co-circles of $H_{\alpha}$
and $H_{\beta}$ provide relations $r^{\prime}_{\alpha}$ and
$r^{\prime}_{\beta}$ in the generators of $\pi_1 (L_3)$.  When the
2-handle dual to $H_{\beta}$ slides over the 2-handle dual to
$H_{\alpha}$, and then back again, $r^{\prime}_{\beta}$ is replaced by
$$r^{\prime}_{\beta}[\mu^{\prime}\overline{\lambda}^{\prime},
r^{\prime}_{\alpha}]^{\overline{\mu}^{\prime}},$$ where
$\lambda^{\prime}$ and $\mu^{\prime}$ describe the homotopy classes of
$\lambda$ and $\mu$ in $\pi_1 (L^3)$.

{\bf 7.} Proposition:  It is possible to choose an arc $\lambda$ which
represents any two given elements in $\pi_1 (L_1)$ and $\pi_1 (L_3)$.
That is, if $j_i : \pi_1 (Q) \rightarrow \pi_1 (L_i)$, $i = 1,3$, then
$$\pi_1 (Q) \stackrel{i_1 \oplus i_3}{\longrightarrow} \pi_1
(L_1) \oplus \pi_1 (L_3)$$ is onto.

Proof:  $\partial Q$ is a collection of tori $T_{\alpha}$;  each
contains loops $\gamma_{\alpha,1}$ and $\gamma_{\alpha,3}$ defined by
$h_{\alpha}(S^1 \times point)$ and $h_{\alpha}(point \times S^1)$ for
the attaching map $h_{\alpha}:S^1 \times B^2 \rightarrow L_1$ for the
2-handle $H_{\alpha}$.  The $\{ \gamma_{\alpha,1} \}$ normally
generate $\pi_1(L_3)$ and represent 0 in $\pi_1(L_1)$, and similarly
the $\{ \gamma_{\alpha,3} \}$ normally generate $\pi_1(L_1)$ and
represent 0 in $\pi_1(L_3)$.  Thus one can represent $(g_1, g_3) \in
\pi_1 (L_1) \oplus \pi_1 (L_3)$ by representing $g_1$ by a loop which
is a product of conjugates of the $\{ \gamma_{\alpha,3} \}$'s, and
similarly $g_3$, and then composing the two loops.

{\bf 8.} Thus, by choosing $\lambda$ and $\mu$ so that they are
homotopic in $L_1$ but are arbitrary in $L_3$, we can slide
$H_{\alpha}$ over $H_{\beta}$ and back so as to replace
$r^{\prime}_{\alpha}$ with $r^{\prime}_{\alpha}$ times any conjugate
of the commutator of any element with $r^{\prime}_{\beta}$, without
changing $r_{\alpha}$.

Recall that we have 2-handles $H_l$, $l = 1, \ldots s$ in $M_{1/2}$
such that the $H_l$, $l = 1, \ldots r$ belong to $A_{1/2}$ and give
relators $r_l$ killing $\pi_1 (L_1)$;  the cocores of the $H_l$, $l =
1, \ldots s$ give relators $r_l^{\prime}$, and the cocores of the
\{$H_{k,i}$\} give relators $r_{k,i}^{\prime}$ which together must
kill $\pi_1 (L_3)$.

Since $H_1 (C_{1/2}) = 0 = \pi_1(C_{1/2})/[\pi_1, \pi_1 ]$, it follows
that the 2-handles in $C_{1/2}$, namely $H_l$, $l = r+1, \ldots s$,
give relators $r_l^{\prime}$ which, modulo $[\pi_1, \pi_1 ]$, kill
$\pi_1 (L_3)$.  More precisely, the relators $r_l^{\prime}$ times a
certain product of conjugates of commutators of arbitrary elements of
$\pi_1 (L_3)$, i.e. $$r_l^{\prime}\prod_j [a_{l,j},
b_{l,j}]^{c_{l,j}}, \;\;\;\; l = r+1, \ldots s,$$ form a set of
elements of $\pi_1 (L_3)$ which normally generate it.  If we had
2-handles $H_{l,j}$ whose cocores represented each of the $b_{l,j}$,
then we could replace $r_l^{\prime}$ by $r_l^{\prime}\prod_j [a_{l,j},
b_{l,j}]^{c_{l,j}}$ by sliding $H_{l,j}$ over $H_l$ and back using
arcs $\lambda_{l,j}$ and $\mu_{l,j}$ where $\mu^{\prime}_{l,j}
\overline{\lambda}^{\prime}_{l,j} = a_{l,j}$ and
$\overline{\mu}^{\prime}_{l,j} = c_{l,j}$, and so that each arc is
trivial in $\pi_1(L_1)$ so that the core of $H_{l,j}$ does not change
its homotopy type.  Having done this replacement, the cocores of the
new $H_l$, $l= r+1, \ldots s$, would kill $\pi_1(L_3)$.

So it suffices to find the 2-handles $H_{l,j}$.  Suppose there are $m$
of the $b_{l,j}$.  Then we introduce $m$ cancelling 1-,2-handle pairs
into the handlebody structure on $M_{1/2}$ and include these pairs in
$A_{1/2}$.  Each $b_{l,j}$ is a product of conjugates of the relators
$r_l^{\prime}$, $ l = 1, \ldots s$, so if we slide the corresponding
2-handles $H_l$ or $H_{k,i}$ over $H_{l,j}$, then the cocore of
$H_{l,j}$ slides over the cocores and ends up representing $b_{l,j}$.
Of course, the core of $H_l$ still kills its original 1-handle, and
sliding $H_{k,i}$ merely changes the isotopy class of the 2-sphere
$S_{k,i}$. (This step is essentially nothing but the observation that
one can always add a cancelling pair of 2-,3-handles where the
2-handle represents any desired word in the 1-handles.)

It may be useful to summarize here the whole construction.  In
$M_{1/2}$, choose a base point, $\ast$, hence a 0-handle, and then $r$
1-handles corresponding to each point of intersection between the
ascending and descending 2-spheres.  Each of these 2-spheres then
provides a 2-handle $H_{k,i}$. Extend this handle structure to
$M_{1/2}$.  Slide 2-handles to get $r$ 2-handles $H_l$ which
homotopically cancel the 1-handles (stablization by cancelling pairs
of 2-, 3-handles to avoid Andrews-Curtis issues may have been
necessary). Add some spare pairs of cancelling 1-,2-handles for later
use.  Inverting $M_{1/2}$ so that 1-handles become 3-handles, etc., we
slide the spare 2-handles over the other 2-handles so that they will
represent certain words, namely the $b_{l,j}$.  Then we slide the
2-handles $H_l^{\prime}$, $l = r+1
\ldots s$, over the spare 2-handles and back so as to create relators
which kill $\pi_1 (L_3)$.

Now $B_{1/2}$ will consist of all of the 1-handles and all their
homotopically cancelling 2-handles, so that $B_{1/2}$ is contractible.
$A_{1/2}$ is $B_{1/2}$ union the $H_{k,i}$'s, as before. Finally
$C_{1/2}$ is the 3-handles union the final version of the $H_l$, $l = 1,
\ldots s$.  Both $A_{1/2}$ and $C_{1/2}$ are simply connected.

We use this new $A_{1/2}$ and proceed to prove Addenda (B), (C) and
(D) as before (it is easy to check in proving (D) that the complement
remains simply connected).  This completes the proof of the Theorem
and all its Addenda.

\bibliography{/u/kirby/Gok/gokova}
\bibliographystyle{plain}

\end{document}